\documentclass{ifacconf}

\usepackage{graphicx}      
\usepackage{natbib}        
\usepackage{comment}  
\usepackage{amsmath,amssymb}
\usepackage{mathtools}
\usepackage{adjustbox}
\usepackage{enumitem}
\usepackage{calc}
\usepackage{booktabs}
\usepackage{multirow}
\usepackage{xcolor}

\DeclarePairedDelimiter\floor{\lfloor}{\rfloor}

\begin{document}
\begin{frontmatter}
\copyright 20XX the authors. This work has been accepted to IFAC for publication under a Creative Commons Licence CC-BY-NC-ND

\title{Hybrid proactive approach for solving maintenance and planning problems in the scenario of Industry 4.0 \thanksref{footnoteinfo}}

\thanks[footnoteinfo]{The authors thank CAPES, CNPq, and FAPEMIG for financial support. M.G.R. also acknowledges partial support of FUNDEP.}

\author[First]{Fernanda F. Alves} 
\author[Second]{Martin G. Ravetti}

\address[First]{Graduate Program in Production Engineering, Federal University of Minas Gerais, Presidente Antonio Carlos Avenue, 6627, 30161-010, Belo Horizonte, Minas Gerais, Brazil (e-mail: fernandafalves@ufmg.br).}
\address[Second]{Department of Computer Science, Federal University of Minas Gerais, Presidente Antonio Carlos Avenue, 6627, Cep 30161-010, Belo Horizonte, Minas Gerais, Brazil (e-mail: gravetti.martin@gmail.com)}

\begin{abstract}                
Behind the concept of Industry 4.0, there are a number of principles and ideas; one of them is the integration of problems of different decision levels. In this work, we integrate maintenance with planning problems, aiming to take full advantage of the production capacity providing immediate delivery of products to customers and avoiding failures. We propose a hybrid approach for solving a maintenance problem integrated with the lot-sizing and scheduling problems. The approach is based on the concepts of robustness and simheuristics, considering preventive, predictive and corrective maintenances in a parallel machine environment. Simulations are performed to consider machine failures. Results indicate that as we increase the robustness parameter at the lot-sizing problem, we obtain lower deviations related to the initial objective function of the lot-sizing problem and lower probabilities of infeasibility after the occurrence of failures. The approach makes it possible to find an average expected result for each of the scenarios analyzed taking into account historical data on the behavior of the failures.
\end{abstract}


\end{frontmatter}

\section{Introduction}

In the context of Industry 4.0, orders will no longer arrive in batches, but individually. Ordering information will be readily available in an online environment. Thus, production planning and scheduling activities should also be carried out quickly. We hypothesize that the lot-sizing model will no longer be considered over a medium-term horizon, but rather over a very short term integrated with scheduling. Also, as systems evolve, decisions will depend less on lot-sizing models, i.e., individual orders will arrive continuously and will be quickly scheduled. Lot-sizing problems can still be used to group products to gain scale and reduce costs. However, as mentioned, lot-sizing decisions will be made much faster. \citet{hofmann2017} states that production planning will be more accurate given that the increased monitoring of flows will lead to a precise demand forecast. Also, as the material flows can be tracked, the utilization of deterministic production planning will be reduced. With real information about the demands, the suppliers will just adjust the production according to it. In this work, we consider an intermediate situation, where it is still necessary to define the production lot sizes, but for this, we consider a short term horizon. In this environment, failures continue to exist, but we have to deal with them much faster.  

We study a lot-sizing and scheduling problem integrated with maintenance activities in a scenario of identical parallel machines. The production process is subject to uncertainties, consisting of machine breakdowns. We propose a hybrid proactive approach for solving the studied problem, which integrates concepts of robustness at the lot-sizing model and concepts of simheuristics to deal with failures at the shop floor. The proposed approach aims to improve the level of customer service, providing products quickly and with quality, see \cite{psarommatis2020}. We did not find works in the literature that dealt with all the problems considered in this work simultaneously. Besides, the proposed resolution method for each specific problem has particularities that are propositions of the present work aiming to consider an emerging yet not clear scenario.

The remainder of this paper is organized as follows: Section \ref{review} presents a literature review, focusing on the integration of maintenance and planning problems. The studied problem is presented in Section \ref{problem}, while in Section \ref{approach} the proposed hybrid approach is explained. Results and test instances are shown in Section \ref{computationalexperiments}, while conclusions are made in Section \ref{conclusion}.

\section{Literature Review}
\label{review}

Most production planning studies do not consider the impact that preventive and corrective maintenances cause in tactical planning. Maintenance activities are usually seen as a disruption and source of costs for the production. However, nowadays, the integration of maintenance and production planning problems is a necessity in the production process (\cite{hafidi2017}). See the work of \cite{ruiz2020} for an example of a maintenance problem in the context of Industry 4.0.

\cite{ettaye2018} study an integrated production planning and maintenance problem aiming to minimize total costs, which consist of production, storage, breaking on the demand, preventive and minimal repair costs. The proposed resolution method determines if we should perform preventive maintenances periodically or perform minimal repairs when unexpected failures occur. Besides, the production quantities are determined considering these interruptions. First, the preventive maintenance problem is solved optimally and then a genetic algorithm (GA) is proposed to solve the production planning based on the results of the first level. Then, the total cost is deduced. This procedure is performed considering all possible periodicities for preventive maintenance in order to find the one resulting in the minimal cost.

In the work of \cite{yalaoui2014}, the authors propose an integrated model for solving a maintenance and production planning problem. The model considers the deterioration of parallel production lines, i.e., the loss in the capacity as time goes on. To improve the computational time, a relaxation technique is proposed for solving moderate instances, while a Fix and Relax heuristic is developed for solving more complex problems, which proved to be efficient when compared to a Lagrangian Relaxation. An optimal strategy for simultaneously solving maintenance and production planning problems is also proposed by \cite{fitouhi2012}. The study considers a single machine subject to random failures, in which preventive maintenance may be performed, while unplanned failures are corrected with minimal repairs. The results show that considering both problems simultaneously reduced the total cost for cyclical and noncyclical preventive replacements. Besides reducing total costs and avoiding the occurrence of failures, the authors state that this synchronization may avoid production delays and re-planning.

\cite{najid2011} propose an integrated production and maintenance planning model considering time windows. The authors state that, usually, in the Enterprise Resource Planning (ERP) systems both problems are considered in separate modules and no integration is performed to minimize production and maintenance costs. They emphasize that it is necessary to consider maintenance as part of the production plan. The proposed model aims to determine the optimal production plan as well as the optimal maintenance periodicity, considering the failure rate function based on the Weibull distribution. To model the problem, the authors consider that the expected failures increase since the last preventive maintenance. They compared the results with a hierarchical approach, in which the production planning model is solved based on the results of the preventive maintenance model. In this latter case, the periodicity of the preventive maintenance is determined only based on expected preventive and corrective costs per period. Results showed that when capacity is moderately loose, the hierarchical approach provides better or equal results compared to the integrated approach. When the capacity is tight, the opposite behavior is observed. When the capacity is too tight, the integrated approach provides better results for all tested instances. However, the instances consider at most 12 products and 24 periods.

Another example of the integration of maintenance and planning problems can be seen in the paper by \cite{ashayeri1996}. Besides, a review of papers considering maintenance integrated with production planning problems can be found in the work of \cite{hafidi2017}, which focuses on the studies dealing with subcontracting constraints. Among the papers cited in this section, the work of \cite{najid2011} is the one that most resembles the present work. However, this paper has significant differences concerning the proposed method and problem.

\section{Problem Description}
\label{problem}

We study a maintenance problem integrated with the lot-sizing and scheduling problems. We consider a parallel machines environment, in which all the jobs can be produced in any machine with the same limited production capacity. The study considers sequence dependent setup times and resumable jobs, i.e., after the unavailable period due to the failure, we may continue the processing of the job which was being produced when the failure occurred. In the non-resumable cases, the job should start all of its processing again. The machines are subject to failures, consisting of machine breakdowns. We consider that failures may occur only during processing and that the failure rate increases over time.

Preemptions are not allowed except in cases where the failure occurs while processing a job and its remaining production quantities have to be finished later. In the algorithms presented in this work, we consider that the production of the same job can occur in the parallel machines at the same time. Such a premise is considered taking into account that a job may be divided into smaller quantities, e.g., quantities of a given weight measure or total amount; but in this case, we must consider the occurrence of sequence dependent setup times on both machines, which may rule out the solution of the best possible solution. In the lot-sizing mathematical model, the job partition in both machines will be less likely to happen, it may occur for example in cases where there is idle production capacity in one machine while the other reached the established capacity limit.

We consider a proactive approach in which we schedule next week's activities. The best solution found is implemented on the production line. If a failure occurs, the line stops, the failure is corrected and the line continues the production based on the original plan. Then, the scheduling is not redone after a failure, we do not go back and redo the planning since we are in a dynamic environment, in which other orders are already being planned. The aim is to have a good notion at the planning level of how the sequence to be executed in the factory will behave if failures occur. We want to send to production a sequence that has a high chance of being feasible, as this sequence will not be revised, i.e., rescheduled as in the reactive approaches. See the paper of \cite{aytug2005} for more information on proactive and reactive approaches.

In a near future, robots will produce customized items quickly on the production line. The tendency is not to have large lots in stock but to produce and immediately deliver products to customers. Then, we believe that inventory costs will be irrelevant on the short-term horizon~\footnote{Unless specific inventory conditions are necessary to stock the products}. Therefore, the objective function to be minimized in the lot-sizing model consists of total backorder quantities, since fast production and delivery will be crucial in this new scenario. Inventory quantities may be used to balance capacity, however, its costs will not be considered. At the scheduling level, our main goal is to minimize the total weighted tardiness.

\section{Hybrid Proactive Approach}
\label{approach}

In this section, we present the proposed hybrid proactive approach. The sets, parameters, and decision variables used in the approach are presented below.

Sets:
\begin{itemize}
	\item $J$: Set of products ($j \in J$)
	\item $M$: Set of machines ($m \in M$)
	\item $S$: Set of subperiods of the planning horizon
	\item $A$: Set of unscheduled products
	\item $Seq_{m}$: Production sequence of machine $m$
\end{itemize}

Parameters:
\begin{itemize}
	\item $K$: Production capacity (in hours)
	\item $\theta_{p}$: Preventive maintenance duration
	\item $\theta_{r}$: Corrective maintenance duration
	\item $D_{j}$: Demand of product $j$
	\item $Sm_{j}$: Average setup time for the production of $j$
	\item $St_{ij}$: Setup time for producing product $j$ after the production of $i$
	\item $p_{j}$: Processing time of product $j$
	\item $d_{j}$: Due date of product $j$
	\item $W_{j}$: Weight of product $j$
	\item $\overline{DF}_{m}$: Average duration of the failures for machine $m$
	\item $SD_{m}$: Standard deviation of the duration of the failures for machine $m$
	\item $P$: Periodicity of the preventive maintenance
	\item $P^{*}$: Optimal periodicity of the preventive maintenance
	\item $L_{m}$: Last product assigned to machine $m$
	\item $CT_{j}$: Completion time of product $j$
	\item $m'$: Machine with the lowest current makespan
	\item $T_{j}$: Tardiness of product $j$
	\item $\beta_{m}$: Shape parameter for machine $m$
	\item $\lambda_{m}$: Scale parameter for machine $m$
	\item $\Delta$: Robustness parameter
\end{itemize}

Decision variables:
\begin{itemize}
	\item $I_{j}$: Continuous variable indicating the quantity stocked of product $j$
	\item $I_{j}^{-}$: Continuous variable indicating the quantity backordered of product $j$
	\item $q_{jm}$: Continuous variable indicating the production quantity of product $j$ on the machine $m$
	\item $w_{jm}$: Binary variable indicating if product $j$ is produced ($w_{jm} = 1$) or not ($w_{jm} = 0$) on machine $m$
\end{itemize}

\subsection{Preventive maintenance problem}

To represent the failure behavior, we use the Weibull distribution. Then, the first step for solving the preventive maintenance problem is to calculate the parameters used in the Weibull distribution, $\beta$ and $\lambda$, which indicates the shape and scale parameters, respectively. For this, based on a set of historical data of time between failures, we use the function $eweibull$ of R package $EnvStats$ to determine the Weibull parameters, which are calculated for each parallel machine $m$ ($\beta_{m}$ and $\lambda_{m}$). Then, we calculate the expected number of failures using Equation \ref{expectedfailures}. In this case, $P$ represents the periodicity of the preventive maintenance and $r_{m}(u)$ represents the failure rate function of the Weibull distribution defined for each machine $m$. 

\begin{equation} \label{expectedfailures}
\int_{0}^{P} r_{m}(u) du = \int_{0}^{P} \left(\dfrac{\beta_{m}}{\lambda_{m}}\right)\left(\dfrac{u}{\lambda_{m}}\right)^{\beta_{m}-1} du  
\end{equation}

In this work, the preventive maintenance problem aims to define the periodicity of the preventive maintenances in order to minimize the total time spent with preventive and corrective maintenances. Then, based on \cite{najid2011}, which considers the total cost, we calculate the total time as shown in the Equation \ref{TotalTime}.

\begin{multline} \label{TotalTime}
Total Time = \sum_{m \in M} \floor*{\dfrac{|S|}{P}}\left(\theta_{p} + \theta_{r}\int_{0}^{P} r_{m}(u) du \right) + \\ + \theta_{r}\int_{0}^{|S|-P\floor*{\dfrac{|S|}{P}}} r_{m}(u) du, \forall{P \in S}
\end{multline}

The second part of Equation \ref{TotalTime} is used only if $\floor*{\dfrac{|S|}{P}}$ is not integer. The total time is calculated for all $P \in S$, selecting for implementation the periodicity $P$ with the lowest total time ($P^{*}$). This output is used as an input for the lot-sizing problem, affecting the production capacity as presented in Equation \ref{Capacity}.

\begin{equation} \label{Capacity}
K = K - \theta_{p}\floor*{\dfrac{|S|}{P}}
\end{equation}

\subsection{Lot-sizing and scheduling problems}
\label{ls}

Based on the results of the preventive maintenance problem, the lot-sizing is solved at optimality using an adaptation for parallel machines of the formulation proposed by \cite{pochetwolsey}. The mathematical model is presented in \eqref{um}-\eqref{sete}. We emphasize that the planning horizon period index is not used in the decision variables because we considered only one planning period. 

\begin{footnotesize}
\begin{flalign}
& \min & & Z = \sum_{j \in J} I_{j}^{-} \label{um} \\
& s.t. & & I_{j} = I_{j0} + \sum_{\mathclap{m \in M}} q_{jm} - D_{j} + I_{j}^{-} - I_{j0}^{-}, & & \forall{j \in J}, \label{dois} \\
& & & \sum_{j \in J} Sm_{j}w_{jm} + \sum_{j \in J} p_{j}q_{jm} \leq K, 					  & & \forall{m \in M}, \label{tres} \\
& & & (p_{j}q_{jm})/K \leq w_{jm} \leq q_{jm},										  	  & & \forall{j \in J}, \forall{m \in M}, \label{quatro} \\
& & & w_{jm} \in \{0,1\}, 																  & & \forall{j \in J}, \forall{m \in M}, \label{cinco}\\
& & & q_{jm} \in \mathbb{R}_{+}, 														  & & \forall{j \in J}, \forall{m \in M}, \label{seis}\\
& & & I_{j}, I_{j}^{-}  \in \mathbb{R}_{+}, 										      & & \forall{j \in J}. \label{sete}
\end{flalign}
\end{footnotesize}

The objective function is presented in Equation \ref{um}, which aims to minimize the total quantity backordered. Constraints \ref{dois} present the conservation flow, with $I_{j0}$ and $I_{j0}^{-}$ indicating the quantities stocked and backordered, respectively, at period $0$. Constraints \ref{tres} present the capacity constraints, in which the sum of total production and setup times must be lower or equal to the capacity of the period. Constraints \ref{quatro} guarantee that average setup times be accounted in the Constraints \ref{tres} only if the product $j$ is produced. Lastly, Constraints \ref{cinco}, \ref{seis} and \ref{sete} present the domain of the decision variables.

At the lot-sizing model, we use concepts of robustness (see \cite{bertsimas2004}), aiming to anticipate the impacts of corrective and predictive maintenances at the tactical level. Then, based on a historical dataset of information about failures, we calculate the average duration of the failures ($\overline{DF}_{m}$) and the standard deviation of the duration of the failures ($SD_{m}$) for each machine $m$. The parameter $\Delta$ is calculated based on $SD_{m}$. We consider three variations to consider robustness for the formulation. In the first case, $\Delta = 0$. In the second case, $\Delta = 0.5SD_{m}$, while in the third case, $\Delta = SD_{m}$. These parameters are added to the Constraints \ref{tres} to generate a robust solution for the problem, as showed in the Constraints \ref{oito}. We expect that, with greater values for $\Delta$, lower probabilities of exceeding the production capacity are found at the shop floor. This is due to the trade-off between the parameter $\Delta$ and the production capacity; as we increase $\Delta$ we use less than the total capacity at the lot-sizing problem, then when we simulate the failures at the scheduling level, the likelihood of exceeding capacity will be lower.   

\begin{equation} \label{oito}
\sum_{j \in J} Sm_{j}w_{jm} + \sum_{j \in J} p_{j}q_{jm} + \overline{DF}_{m} + \Delta \leq K, \forall{m \in M}
\end{equation}

The output of the lot-sizing problem consists of which products will be produced in the period, as well as its respective production quantities, which will represent the inputs for the scheduling problem. As we explained before, in the Industry 4.0 scenario, will be extremely important to deliver the products at their due dates. Then, in the scheduling problem, we aim to determine a production sequence minimizing the total weighted tardiness. Furthermore, the minimization of the makespan is a secondary objective, allowing to produce more products in the period. Here, we propose a heuristic considering job splitting to deal with this problem. Step 1 to Step 5 of the heuristic are based on the \textit{Heuristic 1 (slack-based heuristic)} of \cite{park2012}. However, the procedure of job splitting is entirely different, besides the objective function. Below we present the steps of the modified algorithm.

\begin{enumerate}[labelindent=20pt,labelwidth=\widthof{\ref{last-item}},label=Step \arabic*.,itemindent=0em,leftmargin=!]
  \item Initialize $A$ as all the products to be scheduled, $L_{m} = 0 (\forall{m \in M})$, $CT_{0} = 0$ and $Seq_{m} = \{\emptyset\}$.  
  \item For each product $j \in A$:
    
{\hspace{-0.4cm} a. Set $m'$ to the machine that will be first available (lowest $CT_{L_{m}}$) to process $j$.}
  
{\hspace{-0.4cm} b. Estimate the completion time of product $j$ ($CT_{j}$) when produced on machine $m'$ considering the sequence dependent setup times.}
	\item Calculate the weighted tardiness of product $j$ $T_{j} = W_{j}max(CT_{j} - d_{j},0) \forall{j \in A}$.
	\item Choose $j^{*}$ with the greatest weighted tardiness to be allocated at the first free position of the machine $m'$. Set $L_{m'} = j^{*}$. If more than one product has the same value for $T_{j}$, $j^{*}$ is set as the product $j$ with the greatest weight ($W_{j}$) and lowest due date ($d_{j}$). If more than one product is found, we set $j^{*}$ as the product with the lowest due date. If there is a tie, the tiebreaker is made based on the greatest weight.  
	\item Set $A = A - \{j^{*}\}$ and $Seq_{m'} = Seq_{m'} \cup \{j^{*}\}$. Return to Step 2 until $A = \{\emptyset\}$.
	\item {[Job Splitting]} If one machine exceeds the capacity (infeasible), while the other presents idleness, verify if it is possible to transfer production quantities of the last product of the sequence of the infeasible machine to the idle machine. If it is possible, transfer quantities until the infeasible machine obeys the capacity. After that, balance both machines in a way that both machines present the same ending time, thus minimizing the weighted tardiness of the last product of both machines and the makespan of each machine. Calculate the total weighted tardiness of the solution.
	\item Perform a local search on each machine considering swap neighborhood based on the sequence found before Step 6. Go to Step 6. If the solution found has a weighted total tardiness smaller than the current best solution, save the solution.
	\item Perform a local search on both machines, switching products between them. Go to Step 6. If the solution found has a weighted total tardiness smaller than the current best solution, save the solution.   
\end{enumerate}

OUTPUT: If a solution with the lowest weighted total tardiness and presenting makespans lower or equal to the capacity is found, this solution is an input to the next problem of the proactive approach. However, if only solutions with makespans exceeding the capacity are found, the solution with the lowest weighted total tardiness is the output of the scheduling problem.

The rationale of the scheduling heuristic consists of allocating a product $j$ which results in the greatest weighted tardiness in a machine with the lowest current makespan. When solving the studied scheduling problem with a mathematical model at optimality, only small instances could be solved due to the complexity of the problem. However, the instances solved enabled us to observe a specific characteristic of the solutions that were used for the job splitting step of the heuristic. We observe that, when solving optimally, the job splitting only occurs at the end of the capacity, usually for the last products. This is because we consider sequence dependent setup times. We only produce the same product on both machines if there is not available capacity to produce it in only one machine. We implement this strategy along with local searches, aiming to improve the results obtained by the scheduling heuristic.

\subsection{Predictive and corrective maintenance problem}

We consider that two types of maintenances may occur at the shop-floor beside the preventive maintenance, which consists of the predictive and corrective maintenances. In Industry 4.0, sensors are used to inform when the machine is about to fail by observing its vibration, sound, temperature, among other features. In this case, a signal is sent to the system and predictive maintenance is performed before the failure happens. When the production process already stopped due to the occurrence of a failure, we have to perform corrective maintenance to return the system to its original state. We consider both types of maintenances simulating the moment in which they occur, calculating the expected makespan and the expected weighted total tardiness. These problems are solved based on simheuristics, in which results obtained by the heuristic at the scheduling level are simulated considering the stochastic characteristic of the disruptions. For more information on simheuristics, see the paper of \cite{chica2017}.

\section{Computational Experiments}
\label{computationalexperiments}

The computational tests were run in a computer with Linux operating system, with an Intel\textsuperscript{\textregistered} Core\textsuperscript{TM} i7 processor with 16 GB of RAM. The AMPL R API software and the CPLEX 12.8 solver were used for implementation and resolution. The time limit for the resolution of the lot-sizing model was set to 60 seconds. Preliminary tests showed that within this limit optimal or near-optimal results are found considering the instance set.

\subsection{Instances}
\label{Instances}

For the computational tests, we considered 4, 6, 8, 10, 12, 15, 20, 50, and 100 quantities of products. For each of these quantities, 15 parameter variations were generated considering 3 different variations of the robustness parameter, see Subsection \ref{ls}, totalizing 405 tested instances. We generate the parameters of the problem according to the intervals presented in Table \ref{Distributions}. In the Table, $n$ represents the number of products considered. The demand, due date, and weight are rounded up to the next largest integer, while other parameters are truncated to one decimal place.

\begin{table}[hb]
	\begin{center}
	\caption{Parameter generation.}\label{Distributions}
	\vspace{-0.1cm}
	\begin{tabular}{cccc}
	\hline
	\footnotesize Parameter & \footnotesize Interval \\
	\hline
	Processing time ($p_{j}$) & $U(20/n,20/n+10/n)$ \\
	Setup time ($St_{ij}$) & $U(50/n,50/n+10/n)$ \\
	Due date ($d_{j}$) & $U(max\{j \in J\}(St_{0j} + p_{j}),K)$ \\
	Weight ($W_{j}$) & $U(2,5)$ \\
	Demand ($D_{j}$) & $U(8,12)$ \\
	\hline
	\end{tabular}
	\end{center}
\end{table}

The planning horizon consists of one period of 112 hours. Then, $K = 112$ hours and $|S| = 112$. Two parallel machines are considered. For the simulations, we perform 1000 runs. When generating the predictive and corrective maintenances, the algorithm does not generate numbers in the interval consisting of the moment of preventive maintenance plus a certain interval, stablished as 4 hours, given that the probability that failures occur just after preventive maintenances is low. We spend one hour to perform preventive and predictive maintenances. Lastly, the parameter $\theta_{r}$ was set as equal to $\overline{DF}_{m}$, i.e., the average time to repair the failure based on a historical data set.

\subsection{Results and Discussions}

In this subsection, we present the average results obtained with the proposed approach, see Table \ref{results}. The first column presents the quantities of products. Columns 2, 3, and 4 indicate the average results for Variation 1. In Column 2, we present the average deviation of the expected total weighted tardiness found by the simulation compared to the total weighted tardiness found by the proposed heuristic, which is named $Tard\_Dev$. In Column 3, we present the average deviation of the total backorder quantity after the simulation compared to the initial objective function found by the lot-sizing model before the simulation procedure, named $OF\_Dev$. In Column 4, we present in the first line of the product quantity the average probability of infeasibility after the failure for machine 1, while the second line gives the average probability of infeasibility for machine 2, named $Prob$. This infeasibility occurs if the expected makespan is superior to the production capacity. Columns 5, 6, and 7 presents the results for $Tard\_Dev$, $OF\_Dev$, and $Prob$, respectively, for Variation 2 while columns 8, 9, and 10 present the results for Variation 3. The averages were calculated based on the 15 parameter variations explained in Subsection \ref{Instances}.

\begin{table*}[t]
	\begin{center}
	\caption{Results obtained with the proposed approach.}\label{results}
	\vspace{-0.1cm}
	\resizebox{0.74\textheight}{!}{
	\begin{tabular}{c|ccc|ccc|ccc}
    \toprule
    \multirow{2}[4]{*}{\textbf{Products}} & \multicolumn{3}{c|}{\textbf{Variation 1}} & \multicolumn{3}{c|}{\textbf{Variation 2}} & \multicolumn{3}{c}{\textbf{Variation 3}} \\
\cmidrule{2-10}          & \textbf{Tard\_Dev (\%)} & \textbf{OF\_Dev (\%)} & \textbf{Prob (\%)} & \textbf{Tard\_Dev (\%)} & \textbf{OF\_Dev (\%)} & \textbf{Prob (\%)} & \textbf{Tard\_Dev (\%)} & \textbf{OF\_Dev (\%)} & \multicolumn{1}{c|}{\textbf{Prob (\%)}} \\
    \midrule
    4     & 37.88\% & 6.76\% & 72.54\% & 37.18\% & 5.57\% & 58.94\% & 37.25\% & 4.35\% & 48.79\% \\
          &       &       & 72.35\% &       &       & 64.05\% &       &       & 52.97\% \\
    6     & 22.49\% & 8.46\% & 68.59\% & 27.05\% & 7.58\% & 65.47\% & 30.08\% & 5.83\% & 54.51\% \\
          &       &       & 66.88\% &       &       & 59.12\% &       &       & 57.34\% \\
    8     & 26.78\% & 11.36\% & 79.38\% & 25.09\% & 9.99\% & 72.83\% & 25.87\% & 8.38\% & 66.30\% \\
          &       &       & 81.56\% &       &       & 77.03\% &       &       & 62.39\% \\
    10    & 33.72\% & 11.15\% & 83.89\% & 27.68\% & 9.25\% & 74.40\% & 30.97\% & 7.14\% & 54.95\% \\
          &       &       & 81.81\% &       &       & 68.05\% &       &       & 70.59\% \\
    12    & 29.32\% & 10.54\% & 77.86\% & 26.27\% & 8.61\% & 68.89\% & 24.55\% & 6.63\% & 48.49\% \\
          &       &       & 82.74\% &       &       & 71.65\% &       &       & 66.11\% \\
    15    & 25.04\% & 10.95\% & 82.32\% & 24.55\% & 8.63\% & 70.64\% & 24.77\% & 7.63\% & 62.33\% \\
          &       &       & 81.48\% &       &       & 71.89\% &       &       & 56.19\% \\
    20    & 27.95\% & 10.22\% & 81.14\% & 25.75\% & 8.05\% & 66.88\% & 25.08\% & 6.54\% & 42.52\% \\
          &       &       & 78.69\% &       &       & 69.86\% &       &       & 69.51\% \\
    50    & 39.04\% & 11.91\% & 73.28\% & 39.49\% & 9.94\% & 56.05\% & 42.33\% & 8.43\% & 53.17\% \\
          &       &       & 70.48\% &       &       & 63.87\% &       &       & 46.76\% \\
    100   & 42.26\% & 15.08\% & 70.39\% & 45.50\% & 12.62\% & 63.88\% & 44.50\% & 9.65\% & 48.78\% \\
          &       &       & 69.51\% &       &       & 55.84\% &       &       & 49.29\% \\
    \bottomrule
	\end{tabular}}
	\end{center}
\end{table*}

As we can observe in Table \ref{results}, the average weighted total tardiness obtained after the occurrence of failures is much higher than the results obtained by the scheduling heuristic proposed before the occurrence of failures, ranging from 22.49\% to 42.26\% for Variation 1 for example. However, small differences are observed for the values of $Tard\_Dev$ for most product quantities when different variations in the robustness parameter are considered. Such behavior is because the robustness parameter influence the results obtained by the lot-sizing model. Given that both scheduling heuristics and simulation are performed based on the same results obtained at the upper level by lot-sizing model, it was expected that $Tard\_Dev$ would not be greatly affected by the variation in robustness parameters, which was proven through the results obtained.

On the other hand, when we consider the $OF\_Dev$, we can see in Table \ref{results} that for all product quantities its value decreases as we increase the robustness parameter, i.e., as we move from Variation 1 to Variation 3. This is due to the fact that, when we increase the robustness parameter, we leave a greater gap in the capacity for possible failures. Therefore, in the simulation, such gaps are used, not incurring a large increase in the amount of backorder after the simulations. As we can see, the maximum difference obtained in relation to the original objective function of the lot-sizing model was 15.08\% for 100 products considering Variation 1. The same behavior can be observed for the probabilities of infeasibility after failures ($Prob$). As we increase the robustness parameter, we reduce the probabilities of the expected makespan be greater than the production capacity. The computational time is very low, varying from less than 1 second for 4 products to less than 70 seconds for 100 products.

\section{Conclusions}
\label{conclusion}

In this paper, we proposed a hybrid approach integrating simheuristics and robustness for solving maintenance and production planning problems. The studied scenario consists of a two parallel machines environment, considering sequence dependent setup times. We solve the problems hierarchically, in which first we solve the preventive maintenance problem. Based on the results, the lot-sizing problem is solved, transferring input information for the scheduling problem. Based on this solution, we perform simulations to deal with corrective and predictive maintenances. All problems are studied in the context of the principles and tendencies of Industry 4.0.

Results showed that choosing the more adequate robustness parameter to be used is crucial, given that it directly influences the results of the simulations. We emphasize that the ``more adequate" will depend on the strategic objectives of the company since to guarantee lower probabilities of infeasibility after the failures, we need to be willing to reduce the use of production capacity at the lot-sizing level. For future works, we aim to consider machine learning tools to deal with the predictive maintenance.

\bibliography{References} 

\end{document}